\documentclass{manuscript}

\begin{document}
\title{Normalized solutions for logarithmic $p$-Laplacian equations on graphs}
\author{%
    Wenyuan Ma$^{1}$\thanks{Corresponding author. E-mail: \href{mailto:wenyuanma1203@ruc.edu.cn}{\nolinkurl{wenyuanma1203@ruc.edu.cn}}}, Anfeng Yuan$^{2}$\\
    \affiliation{$^{1}$Department of Mathematics, Renmin University of China, Beijing 100872, China}\\
    \affiliation{$^{2}$Institute of Fundamental and Interdisciplinary Sciences, Beijing Union University,}\\
    \affiliation{Beijing 100101, China}
}
\maketitle

%% ── Abstract ─────────────────────────────────────────
% Abstract
{\noindent\small{\bf Abstract }
    We study a logarithmic $p$-Laplacian equation with prescribed $L^2$ mass on connected locally finite weighted graphs, with $p > 2$.
    By minimizing the energy, we obtain a positive solution for every prescribed mass $m > 0$ on finite graphs and, under additional assumptions on the vertex measure and the potential, on infinite graphs.
    We describe the behavior of positive solutions and their Lagrange multipliers as $m \to 0^+$ or $m \to +\infty$.
    In particular, as $m \to 0^+$, the minimizers rescaled to unit mass converge along subsequences to functions supported at a single vertex.
}

\vspace{1ex}

\vspace{1ex}

%% ── Sections ─────────────────────────────────────────
\section{Introduction}
\label{sec:introduction}
Carles and Gallagher~\cite{Carles2018} note that logarithmic Schrödinger equations have been used as models in nonlinear wave mechanics and nonlinear optics. For the defocusing equation with nonzero $H^1$ initial data of finite second moment, they prove weak $L^1$ convergence of the rescaled density to a universal Gaussian as time tends to infinity. For the standing wave equation
\begin{equation*}
    -\Delta u + \omega u = u \log u^2, \qquad x \in \mathbb{R}^N,
\end{equation*}
where $\omega \in \mathbb{R}$ and $N \geq 3$, the energy can be infinite on $H^1(\mathbb{R}^N)$. Using nonsmooth critical point theory, d'Avenia, Montefusco and Squassina~\cite{DAvenia2014} obtain infinitely many radial weak solutions. For continuous coercive potentials, Ji and Szulkin~\cite{Ji2016} obtain infinitely many solutions. For such potentials, Shuai~\cite{Shuai2019} obtains positive and sign-changing solutions using directional derivatives and constrained minimization.

Grigor'yan, Lin and Yang use variational methods to study Yamabe type equations on bounded domains of locally finite graphs and on finite graphs~\cite{Grigoryan2016c}. They also give conditions for semilinear equations on locally finite graphs to admit positive solutions~\cite{GrigoryanLinYang2017}. Zhang and Zhao~\cite{Zhang2018b} prove subsequential convergence of ground states to a Dirichlet ground state on the potential well as it deepens. Under their respective hypotheses for graph $p$-Laplacian equations, Han and Shao~\cite{Han2021a} obtain positive solutions and ground states for $p \geq 2$, and Shao, Yang and Zhao~\cite{Shao2026a} obtain positive ground states for $p > 1$.

Chang, Wang and Yan~\cite{ChangWangYan2023} obtain ground states of
\begin{equation*}
    -\Delta u + a(x)u = u \log u^2, \qquad x \in V,
\end{equation*}
on locally finite graphs under conditions allowing $a$ to change sign. He and Ji~\cite{HeJi2024} obtain infinitely many solutions on lattice graphs for coercive potentials and ground states for periodic potentials. In Section~6 they also study a logarithmic $p$-Laplacian equation without a mass constraint. Chang et al.~\cite{Chang2023a} study least energy sign-changing solutions on locally finite graphs for sufficiently deep potential wells. For logarithmic Dirichlet problems on finite domains, Shao, Yang and Zhao~\cite{Shao2024c} study ground states, sign-changing solutions and the limit as the nonlinear exponent decreases to $2$.

For prescribed mass in $\mathbb{R}^N$, Alves and Ji~\cite{Alves2024} study multiplicity for a logarithmic Schrödinger equation. For sufficiently small $L^2$ mass, Wang et al.~\cite{Wang2021} obtain normalized solutions of a $p$-Laplacian equation with $L^2$-supercritical power growth. Shen and Squassina~\cite{Shen2025} study existence and concentration for logarithmic $p$-Laplacian equations with prescribed $L^p$ norm. On lattice graphs, He, Ji and Tao~\cite{He2025} study global minimizers at fixed mass.

On connected finite graphs, Yang and Zhao~\cite{YangZhao2024} prove that for every $q > 2$ and $m > 0$ there is a positive solution $u_m$ of
\begin{equation*}
    -\Delta u + \lambda_m u = u^{q - 1}, \qquad \int_V u^2 \, \mathrm{d}\mu = m.
\end{equation*}
On infinite locally finite graphs they obtain positive solutions of
\begin{equation*}
    -\Delta u + \lambda_m h(x)u = u^{q - 1}, \qquad \int_V h(x)u^2 \, \mathrm{d}\mu = m,
\end{equation*}
under their assumptions (c1)--(c3). For a fixed graph, positive weight $h$ and vertex in (c3), that condition bounds $m$ away from zero. They study $u_m/\sqrt{m}$ in both settings. On infinite graphs they further assume $\lambda_m \geq -c_0$ as $m \to 0^+$ and $\int_V |\nabla u_m|^2 \, \mathrm{d}\mu \leq c_\infty m$ as $m \to +\infty$, with $c_0 \geq 0$ and $c_\infty > 0$ independent of $m$~\cite[Theorem~1.4]{YangZhao2024}.

Let $G = (V, E)$ be a connected, locally finite weighted graph with vertex set $V$, edge set $E$ and a positive vertex measure $\mu$. We consider the following normalized logarithmic $p$-Laplacian problem
\begin{equation}
    \begin{cases}
        -\Delta_p u + [a(x) + \lambda_m] u = b u \log u, & x \in V, \\[1mm]
        u > 0,                                           & x \in V, \\[1mm]
        \displaystyle\int_V u^2 \, \mathrm{d}\mu = m,
    \end{cases}
    \label{eq:normalized-logarithmic-p-laplacian}
\end{equation}
where $a \colon V \to \mathbb{R}$ and $b > 0$ are prescribed, and we assume $p > 2$ and $m > 0$. The unknowns are a positive function $u$ on $V$ and a real number $\lambda_m$, which is the Lagrange multiplier associated with the mass constraint $\int_V u^2 \, \mathrm{d}\mu = m$.

We prove the existence of positive normalized solutions and study the behavior of these solutions and their Lagrange multipliers as $m \to 0^+$ and $m \to +\infty$.

To state our main results, we first introduce some notation. For each edge $xy \in E$, let $\omega_{xy} = \omega_{yx} > 0$ be its weight. We write $y \sim x$ if $xy \in E$.

For a function $u \colon V \to \mathbb{R}$, its gradient at $x \in V$ is defined by
\begin{equation*}
    \nabla u(x) \coloneqq \left(\sqrt{\frac{\omega_{xy_i}}{2\mu(x)}} (u(y_i) - u(x))\right)_{i = 1}^{\ell_x},
\end{equation*}
where $y_1, \dotsc, y_{\ell_x}$ are the neighbors of $x$ and $\ell_x$ is their number. The length of the gradient is given by
\begin{equation*}
    |\nabla u|(x) = \left[\frac{1}{2\mu(x)} \sum_{y \sim x} \omega_{xy} (u(y) - u(x))^2\right]^{1/2}.
\end{equation*}
For two functions $u, v \colon V \to \mathbb{R}$, the corresponding gradient form is given by
\begin{equation*}
    \nabla u(x) \cdot \nabla v(x) = \frac{1}{2\mu(x)} \sum_{y \sim x} \omega_{xy} (u(y) - u(x)) (v(y) - v(x)).
\end{equation*}
For $p > 2$, the $p$-Laplacian of $u$ at $x \in V$ is defined by
\begin{equation*}
    \Delta_p u(x) \coloneqq \frac{1}{2\mu(x)} \sum_{y \sim x} \omega_{xy} \left(|\nabla u|(x)^{p - 2} + |\nabla u|(y)^{p - 2}\right) (u(y) - u(x)).
\end{equation*}
Integration with respect to $\mu$ is defined by
\begin{equation*}
    \int_V f \, \mathrm{d}\mu \coloneqq \sum_{x \in V} \mu(x) f(x),
\end{equation*}
whenever the sum is absolutely convergent.

For $1 \leq q < +\infty$, we define
\begin{equation*}
    \|u\|_{L^q(V)} \coloneqq \left(\int_V |u|^q \, \mathrm{d}\mu\right)^{1/q}, \qquad \|u\|_{L^\infty(V)} \coloneqq \sup_{x \in V} |u(x)|.
\end{equation*}
For $1 \leq q \leq +\infty$, we set $L^q(V) \coloneqq \{u \colon V \to \mathbb{R} \colon \|u\|_{L^q(V)} < +\infty\}$.

Throughout this paper, $s \log s$, $s^2 \log s$ for $s > 0$ and $s \log |s|$, $s^2 \log |s|$ for $s \neq 0$, are continuously extended to $s = 0$ by assigning them the value zero.

For a function $u \colon V \to \mathbb{R}$ for which the following integrals are absolutely convergent, we define the energy functional
\begin{equation}
    J(u) \coloneqq \frac{1}{p} \int_V |\nabla u|^p \, \mathrm{d}\mu + \frac{1}{2} \int_V a(x) u^2 \, \mathrm{d}\mu - \frac{b}{2} \int_V u^2 \log |u| \, \mathrm{d}\mu + \frac{b}{4} \int_V u^2 \, \mathrm{d}\mu.
    \label{eq:energy-functional}
\end{equation}

We first consider finite graphs. In this case, we set $\mathcal{H}_0 \coloneqq \{u \colon V \to \mathbb{R}\}$ and equip it with the norm
\begin{equation*}
    \|u\|_{\mathcal{H}_0} \coloneqq \left(\int_V |\nabla u|^p \, \mathrm{d}\mu\right)^{1/p} + \left(\int_V |u|^2 \, \mathrm{d}\mu\right)^{1/2}.
\end{equation*}
In this case, $J$ takes finite values on $\mathcal{H}_0$. Then we have the following existence result.

\begin{theorem}
    \label{thm:finite-graph-existence}
    For every $m > 0$, there exist a function $u_m \in \mathcal{H}_0$ and a corresponding Lagrange multiplier $\lambda_m \in \mathbb{R}$ satisfying \eqref{eq:normalized-logarithmic-p-laplacian}.
\end{theorem}

The next result describes the behavior of these solutions as the mass tends to zero or infinity.
\begin{theorem}
    \label{thm:finite-graph-mass-limits}
    Under the assumptions of Theorem~\ref{thm:finite-graph-existence}, let $u_m$ and its Lagrange multiplier $\lambda_m$ satisfy \eqref{eq:normalized-logarithmic-p-laplacian} for each $m > 0$.
    \begin{enumerate}[(i)]
        \item As $m \to 0^+$, up to a subsequence, we have
        \begin{equation*}
            \frac{u_m}{\sqrt{m}} \to v \quad \text{in } \mathcal{H}_0, \qquad \lambda_m - \frac{b}{2} \log m \to \theta,
        \end{equation*}
        where, for some nonempty set $K \subseteq V$,
        \begin{equation}
            v(x) = \frac{\mathrm{e}^{a(x)/b} \mathbf{1}_K(x)}{\sqrt{\sigma_K}}, \qquad \theta = -\frac{b}{2} \log \sigma_K, \qquad \sigma_K \coloneqq \sum_{y \in K} \mu(y) \mathrm{e}^{2a(y)/b}.
            \label{eq:small-mass-explicit-limit}
        \end{equation}

        Moreover, if each $u_m$ is a minimizer obtained in the proof of Theorem~\ref{thm:finite-graph-existence}, then
        \begin{equation}
            v(x) = \begin{cases}
                \mu(x_*)^{-1/2}, & x = x_*,    \\[1mm]
                0,               & x \neq x_*,
            \end{cases}
            \label{eq:small-mass-minimizer-limit}
        \end{equation}
        where $x_* \in V$ satisfies $\mu(x_*) \mathrm{e}^{2a(x_*)/b} = \min\limits_{x \in V} \mu(x) \mathrm{e}^{2a(x)/b}$.

        \item As $m \to +\infty$, we have, for any solution,
        \begin{equation*}
            \frac{u_m}{\sqrt{m}} \to \mu(V)^{-1/2} \quad \text{in } \mathcal{H}_0
        \end{equation*}
        and
        \begin{equation*}
            \lambda_m - \frac{b}{2} \log m \to -\frac{1}{\mu(V)} \int_V a(x) \, \mathrm{d}\mu - \frac{b}{2} \log \mu(V),
        \end{equation*}
        where $\mu(V) \coloneqq \int_V 1 \, \mathrm{d}\mu$.
    \end{enumerate}
\end{theorem}

We now consider infinite graphs. In addition to the graph and parameter conditions above, we impose the following assumptions.
\begin{enumerate}[({A}1)]
    \item The measure satisfies $\mu(x) \geq \mu_0 > 0$ for every $x \in V$.
    \item For every $R \in \mathbb{R}$, the set
    \begin{equation*}
        \{x \in V \colon a(x) \leq R\}
    \end{equation*}
    is finite.
\end{enumerate}

Fix $a_0 > 0$. By (A2), the set $\{x \in V \colon a(x) \leq a_0\}$ is finite, so there exists $c \in \mathbb{R}$ such that $a(x) \geq c$ for every $x \in V$. Choose $C \geq \max\{0, a_0 - c\}$ and set
\begin{equation*}
    \widetilde a(x) \coloneqq a(x) + C \geq a_0, \qquad \widetilde\lambda_m \coloneqq \lambda_m - C.
\end{equation*}
Replacing $a$ and $\lambda_m$ by $\widetilde a$ and $\widetilde\lambda_m$ leaves problem \eqref{eq:normalized-logarithmic-p-laplacian} unchanged. For convenience, we additionally assume
\begin{enumerate}[({A}1)]
    \setcounter{enumi}{2}
    \item The potential satisfies $a(x) \geq a_0 > 0$ for every $x \in V$.
\end{enumerate}

Let $C_c(V) \coloneqq \{u \colon V \to \mathbb{R} \colon \operatorname{supp} u \text{ is finite}\}$. For infinite graphs, we define $\mathcal{H}$ as the completion of $C_c(V)$ under the norm
\begin{equation*}
    \|u\|_{\mathcal{H}} \coloneqq \left(\int_V |\nabla u|^p \, \mathrm{d}\mu\right)^{1/p} + \left(\int_V a(x) |u|^2 \, \mathrm{d}\mu\right)^{1/2},
\end{equation*}
that is, $\mathcal{H} \coloneqq \overline{C_c(V)}^{\|\cdot\|_{\mathcal{H}}}$.
The same formula \eqref{eq:energy-functional} defines $J \colon \mathcal{H} \to \mathbb{R} \cup \{+\infty\}$, where $J(u) \coloneqq +\infty$ whenever the logarithmic integral is not absolutely convergent.

We also define $\mathcal{W}$ as the completion of $C_c(V)$ under the norm
\begin{equation*}
    \|u\|_{\mathcal{W}} \coloneqq \left(\int_V |\nabla u|^p \, \mathrm{d}\mu\right)^{1/p} + \left(\int_V |u|^2 \, \mathrm{d}\mu\right)^{1/2},
\end{equation*}
that is, $\mathcal{W} \coloneqq \overline{C_c(V)}^{\|\cdot\|_{\mathcal{W}}}$.
Then $\mathcal{H} \hookrightarrow \mathcal{W}$.

For infinite graphs, a solution of problem \eqref{eq:normalized-logarithmic-p-laplacian} is a pair $(u, \lambda_m) \in \mathcal{H} \times \mathbb{R}$ satisfying the problem pointwise, which is equivalent to
\begin{equation}
    \int_V |\nabla u|^{p - 2} \nabla u \cdot \nabla \phi \, \mathrm{d}\mu + \int_V [a(x) + \lambda_m] u \phi \, \mathrm{d}\mu = b \int_V u \phi \log u \, \mathrm{d}\mu
    \label{eq:infinite-graph-variational-identity}
\end{equation}
for every $\phi \in C_c(V)$.

\begin{theorem}
    \label{thm:infinite-graph-existence}
    Assume (A1)--(A3). For every $m > 0$, there exist a function $u_m \in \mathcal{H}$ and a corresponding Lagrange multiplier $\lambda_m \in \mathbb{R}$ satisfying \eqref{eq:normalized-logarithmic-p-laplacian}.
\end{theorem}

For the large-mass limit, set
\begin{equation}
    \kappa \coloneqq \inf_{\substack{w \in C_c(V) \\ \int_V w^2 \, \mathrm{d}\mu = 1}} \int_V |\nabla w|^p \, \mathrm{d}\mu.
    \label{eq:infinite-graph-gradient-infimum}
\end{equation}
The following theorem describes both mass limits.
\begin{theorem}
    \label{thm:infinite-graph-mass-limits}
    Assume (A1)--(A3), and let $u_m$ and its Lagrange multiplier $\lambda_m$ satisfy \eqref{eq:normalized-logarithmic-p-laplacian} for each $m > 0$.
    \begin{enumerate}[(i)]
        \item As $m \to 0^+$, up to a subsequence, either
        \begin{equation*}
            \lambda_m - \frac{b}{2} \log m \to \theta \in \mathbb{R}, \qquad \frac{u_m}{\sqrt{m}} \to v \quad \text{in } L^2(V),
        \end{equation*}
        where $v$ and $\theta$ are given by \eqref{eq:small-mass-explicit-limit} for some nonempty finite set $K \subseteq V$, or
        \begin{equation*}
            \lambda_m - \frac{b}{2} \log m \to -\infty, \qquad \frac{u_m(x)}{\sqrt{m}} \to 0 \quad \text{for every } x \in V.
        \end{equation*}

        If each $u_m$ is a minimizer obtained in the proof of Theorem~\ref{thm:infinite-graph-existence}, then only the first alternative occurs and $v$ is given by \eqref{eq:small-mass-minimizer-limit}, where $x_* \in V$ minimizes $\mu(x) \mathrm{e}^{2a(x)/b}$.

        \item As $m \to +\infty$, up to a subsequence, we have
        \begin{equation*}
            \frac{u_m(x)}{\sqrt{m}} \to v(x) \quad \text{for every } x \in V,
        \end{equation*}
        where $v \in L^2(V)$. Either $v \equiv 0$, or $v$ satisfies
        \begin{equation}
            \begin{cases}
                -\Delta_p v - \kappa v = 0, & x \in V, \\[1mm]
                v > 0, & x \in V, \\[1mm]
                \displaystyle\int_V v^2 \, \mathrm{d}\mu = 1,
            \end{cases}
            \label{eq:infinite-graph-large-mass-limit}
        \end{equation}
        where $\kappa > 0$ is defined by \eqref{eq:infinite-graph-gradient-infimum}. In this case, along the same subsequence,
        \begin{equation*}
            \frac{u_m}{\sqrt{m}} \to v \quad \text{in } L^2(V), \qquad \frac{1}{m^{p/2 - 1}} \left(\lambda_m - \frac{b}{2} \log m\right) \to -\kappa.
        \end{equation*}

        Moreover, if each $u_m$ is a minimizer obtained in the proof of Theorem~\ref{thm:infinite-graph-existence} and $\kappa > 0$, then $v \not\equiv 0$.
    \end{enumerate}
\end{theorem}

The paper is organized as follows. In Section~\ref{sec:finite-graph}, we prove Theorems~\ref{thm:finite-graph-existence} and~\ref{thm:finite-graph-mass-limits} for finite graphs. Section~\ref{sec:infinite-graph} contains the proofs of Theorems~\ref{thm:infinite-graph-existence} and~\ref{thm:infinite-graph-mass-limits} for infinite graphs.

\section{The finite graph case}
\label{sec:finite-graph}

Throughout this section, $G = (V, E)$ is a connected finite graph satisfying the assumptions in Theorem~\ref{thm:finite-graph-existence}. We prove the theorem by minimizing the energy under the mass constraint.

We introduce the functional $I \colon \mathcal{H}_0 \to \mathbb{R}$ by
\begin{equation*}
    I(u) \coloneqq \frac{1}{2} \int_V a(x) u^2 \, \mathrm{d}\mu - \frac{b}{2} \int_V u^2 \log |u| \, \mathrm{d}\mu.
\end{equation*}
We have $s^2 \log |s| \in C^1(\mathbb{R})$, with derivative $2s \log |s| + s$. We shall minimize the energy functional $J$ defined in \eqref{eq:energy-functional} on the set
\begin{equation*}
    \mathcal{S}_m^0 \coloneqq \left\{u \in \mathcal{H}_0 \colon \int_V u^2 \, \mathrm{d}\mu = m\right\}.
\end{equation*}

We now prove Theorem~\ref{thm:finite-graph-existence}, following the argument of Yang and Zhao \cite{YangZhao2024}.

\begin{proof}[Proof of Theorem~\ref{thm:finite-graph-existence}]
    Write $V = \{x_1, \dotsc, x_k\}$, where $k \in \mathbb{N}_+$. Define $\mathbf{y} \coloneqq (u(x_1), \dotsc, u(x_k))$, then $\mathcal{H}_0$ can be identified with $\mathbb{R}^k$. Hence, $u \in \mathcal{S}_m^0$ if and only if
    \begin{equation*}
        \mathbf{y} \in A_m \coloneqq \left\{(y_1, \dotsc, y_k) \in \mathbb{R}^k \colon \sum_{i = 1}^k \mu(x_i) y_i^2 = m\right\}.
    \end{equation*}
    Define $F \colon \mathbb{R}^k \to \mathbb{R}$ by $F(\mathbf{y}) \coloneqq J(u)$, where $u(x_i) \coloneqq y_i$ for $i = 1, \dotsc, k$. And we have $F \in C^1(\mathbb{R}^k, \mathbb{R})$.

    Since $(\sqrt{m/\mu(x_1)}, 0, \dotsc, 0) \in A_m$, the set $A_m$ is nonempty. The map $\mathbf{y} \mapsto \sum\limits_{i = 1}^k \mu(x_i) y_i^2$ is continuous, so $A_m$ is closed. Moreover, every $\mathbf{y} \in A_m$ satisfies
    \begin{equation*}
        \|\mathbf{y}\|_{\mathbb{R}^k}^2 = \sum_{i = 1}^k y_i^2 \leq \frac{m}{\displaystyle\min_{1 \leq i \leq k} \mu(x_i)}.
    \end{equation*}
    Thus $A_m$ is bounded and hence compact. Consequently, $F$ attains its minimum at some $\mathbf{y}^m = (y_1^m, \dotsc, y_k^m) \in A_m$. Setting $u_m(x_i) \coloneqq y_i^m$, we obtain $u_m \in \mathcal{S}_m^0$ and
    \begin{equation*}
        J(u_m) = \inf_{u \in \mathcal{S}_m^0} J(u).
    \end{equation*}

    Next, we show that the minimizer may be chosen nonnegative. For every $u \in \mathcal{H}_0$ and $x \in V$, we have
    \begin{equation*}
        \begin{aligned}
            |\nabla |u||^2(x) & = \frac{1}{2\mu(x)} \sum_{y \sim x} \omega_{xy} (|u(y)| - |u(x)|)^2                   \\
                              & \leq \frac{1}{2\mu(x)} \sum_{y \sim x} \omega_{xy} (u(y) - u(x))^2 = |\nabla u|^2(x).
        \end{aligned}
    \end{equation*}
    The remaining terms in $J$ depend only on $|u|$, and therefore $J(|u|) \leq J(u)$. Since $|u_m| \in \mathcal{S}_m^0$, we have
    \begin{equation*}
        J(|u_m|) = \inf_{u \in \mathcal{S}_m^0} J(u).
    \end{equation*}
    Replacing $u_m$ by $|u_m|$ if necessary, we may assume that $u_m \geq 0$ on $V$.

    We now derive the Euler--Lagrange equation. Define
    \begin{equation*}
        Q(u) \coloneqq \frac{1}{2} \int_V u^2 \, \mathrm{d}\mu, \qquad \mathcal{S}_m^0 = \left\{u \in \mathcal{H}_0 \colon Q(u) = \frac{m}{2}\right\}.
    \end{equation*}
    Since $\langle Q'(u_m), u_m\rangle = m > 0$, we have $Q'(u_m) \neq 0$. The Lagrange multiplier method therefore gives $\lambda_m \in \mathbb{R}$ such that
    \begin{equation*}
        \langle J'(u_m), \phi\rangle + \lambda_m \langle Q'(u_m), \phi\rangle = 0 \qquad \text{for every } \phi \in \mathcal{H}_0,
    \end{equation*}
    where $\langle\cdot, \cdot\rangle$ denotes the duality pairing between $\mathcal{H}_0^*$ and $\mathcal{H}_0$. That is,
    \begin{equation}
        \int_V |\nabla u_m|^{p - 2} \nabla u_m \cdot \nabla \phi \, \mathrm{d}\mu + \int_V [a(x) + \lambda_m] u_m \phi \, \mathrm{d}\mu = b \int_V u_m \phi \log u_m \, \mathrm{d}\mu.
        \label{eq:finite-graph-variational-equation}
    \end{equation}
    Integration by parts and taking $\phi$ to be the indicator function of $\{x\}$ yield
    \begin{equation}
        -\Delta_p u_m(x) + [a(x) + \lambda_m] u_m(x) = b u_m(x) \log u_m(x), \qquad x \in V.
        \label{eq:finite-graph-euler-lagrange}
    \end{equation}

    We claim that $u_m(x) > 0$ for every $x \in V$. Suppose, for contradiction, that $u_m(x_0) = 0$ for some $x_0 \in V$. By \eqref{eq:finite-graph-euler-lagrange}, $\Delta_p u_m(x_0) = 0$. Thus
    \begin{equation*}
        \frac{1}{2\mu(x_0)} \sum_{y \sim x_0} \omega_{x_0 y} \left(|\nabla u_m|(x_0)^{p - 2} + |\nabla u_m|(y)^{p - 2}\right) u_m(y) = 0.
    \end{equation*}
    If $u_m(y) > 0$ for some $y \sim x_0$, then $|\nabla u_m|(x_0) > 0$, and we have $\Delta_p u_m(x_0) > 0$. This is impossible, so $u_m(y) = 0$ for every neighbor of $x_0$. Repeating this argument along paths and using the connectedness of $G$, we obtain $u_m \equiv 0$ on $V$, contradicting $\int_V u_m^2 \, \mathrm{d}\mu = m > 0$. Hence $u_m$ is strictly positive on $V$.

    Finally, taking $\phi = u_m$ in \eqref{eq:finite-graph-variational-equation}, we obtain
    \begin{equation*}
        \lambda_m = -\frac{1}{m} \left(\int_V |\nabla u_m|^p \, \mathrm{d}\mu + 2I(u_m)\right),
    \end{equation*}
    which is finite because $V$ is finite. Together with $u_m \in \mathcal{S}_m^0$ and $u_m > 0$ on $V$, \eqref{eq:finite-graph-euler-lagrange} proves that $(u_m, \lambda_m)$ solves \eqref{eq:normalized-logarithmic-p-laplacian}.
\end{proof}

We next prove Theorem~\ref{thm:finite-graph-mass-limits}, following the approach used in the proof of \cite[Theorem~1.2]{YangZhao2024}.

\begin{proof}[Theorem~\ref{thm:finite-graph-mass-limits}(i)]
    Let $\{m_n\}$ be an arbitrary sequence of positive masses such that $m_n \to 0^+$ as $n \to \infty$, $n \in \mathbb{N}_+$. Set $v_n \coloneqq u_{m_n}/\sqrt{m_n}$. The constraint gives
    \begin{equation*}
        v_n(x) > 0 \quad \text{for every } x \in V, \qquad \int_V v_n^2 \, \mathrm{d}\mu = 1.
    \end{equation*}
    In particular, $v_n \in \mathcal{S}_1^0$ and $0 < v_n(x) \leq \mu(x)^{-1/2}$ for every $x \in V$. Since $V$ is finite, $\mathcal{S}_1^0$ is compact in $\mathcal{H}_0$. Thus, after passing to a subsequence, there exists $v \in \mathcal{H}_0$ such that $v_n \to v$ in $\mathcal{H}_0$ and uniformly on $V$. Consequently,
    \begin{equation*}
        v \geq 0 \quad \text{on } V, \qquad \int_V v^2 \, \mathrm{d}\mu = 1.
    \end{equation*}

    By the homogeneity of the $p$-Laplacian, $\Delta_p(\sqrt{m_n} v_n) = m_n^{(p - 1)/2} \Delta_p v_n$. Substituting $u_{m_n} = \sqrt{m_n} v_n$ into \eqref{eq:normalized-logarithmic-p-laplacian}, we obtain
    \begin{equation*}
        -m_n^{p/2 - 1} \Delta_p v_n + [a(x) + \lambda_{m_n}] v_n = b v_n \log v_n + \frac{b}{2} v_n \log m_n, \qquad x \in V.
    \end{equation*}
    Define
    \begin{equation*}
        \tau_n \coloneqq m_n^{p/2 - 1}, \qquad \theta_n \coloneqq \lambda_{m_n} - \frac{b}{2} \log m_n.
    \end{equation*}
    Then the equation becomes
    \begin{equation}
        -\tau_n \Delta_p v_n + [a(x) + \theta_n] v_n = b v_n \log v_n, \qquad x \in V.
        \label{eq:finite-graph-mass-rescaled}
    \end{equation}

    Multiplying \eqref{eq:finite-graph-mass-rescaled} by $v_n$ and integrating over $V$, we find
    \begin{equation}
        \theta_n = -2I(v_n) - \tau_n \int_V |\nabla v_n|^p \, \mathrm{d}\mu.
        \label{eq:finite-graph-mass-multiplier}
    \end{equation}
    Since $V$ is finite, the uniform convergence $v_n \to v$ yields
    \begin{equation*}
        I(v_n) \to I(v).
    \end{equation*}
    The definition of the gradient also gives
    \begin{equation*}
        \int_V |\nabla v_n|^p \, \mathrm{d}\mu \to \int_V |\nabla v|^p \, \mathrm{d}\mu.
    \end{equation*}
    In particular, these gradient integrals are bounded. Since $p > 2$, we have $\tau_n \to 0$, and \eqref{eq:finite-graph-mass-multiplier} implies
    \begin{equation*}
        \theta_n \to \theta \coloneqq -2I(v) \in \mathbb{R}.
    \end{equation*}

    Since $\tau_n \to 0$ and $\Delta_p v_n(x) \to \Delta_p v(x)$ at each vertex, passing to the limit in \eqref{eq:finite-graph-mass-rescaled} gives
    \begin{equation*}
        [a(x) + \theta] v(x) = b v(x) \log v(x), \qquad x \in V.
    \end{equation*}
    Let $K \coloneqq \operatorname{supp} v$, which is nonempty since $\int_V v^2 \, \mathrm{d}\mu = 1$. On $K$, dividing by $v(x) > 0$ gives $v(x) = \mathrm{e}^{(a(x) + \theta)/b}$, and the mass constraint yields
    \begin{equation}
        \sum_{x \in K} \mu(x) v(x)^2 = \mathrm{e}^{2\theta/b} \sigma_K = 1.
        \label{eq:small-mass-normalization}
    \end{equation}
    Solving for $\theta$ and substituting into $v(x)$ gives \eqref{eq:small-mass-explicit-limit}. Since $\theta = -2I(v)$, we also have
    \begin{equation}
        I(v) = \frac{b}{4} \log \sigma_K.
        \label{eq:finite-graph-limit-energy}
    \end{equation}

    Suppose now that each $u_m$ minimizes $J$ on $\mathcal{S}_m^0$. We next show that the subsequential limit $v$ minimizes $I$ on $\mathcal{S}_1^0$. 
    
    For any $w \in \mathcal{S}_1^0$, set $u \coloneqq \sqrt{m} w$. Using $\int_V w^2 \, \mathrm{d}\mu = 1$, we obtain
    \begin{equation}
        J(u) = J(\sqrt{m} w) = \frac{m^{p/2}}{p} \int_V |\nabla w|^p \, \mathrm{d}\mu + m I(w) + \frac{bm}{4} (1 - \log m).
        \label{eq:finite-graph-scaled-energy}
    \end{equation}

    By the minimality of $u_{m_n}$, we have
    \begin{equation*}
        J(u_{m_n}) = J(\sqrt{m_n} v_n) \leq J(\sqrt{m_n} w).
    \end{equation*}
    Consequently,
    \begin{equation}
        I(v_n) + \frac{\tau_n}{p} \int_V |\nabla v_n|^p \, \mathrm{d}\mu \leq I(w) + \frac{\tau_n}{p} \int_V |\nabla w|^p \, \mathrm{d}\mu.
        \label{eq:finite-graph-minimizer-comparison}
    \end{equation}
    Since $V$ is finite, letting $n \to \infty$ yields
    \begin{equation*}
        I(v) \leq I(w).
    \end{equation*}
    As $w \in \mathcal{S}_1^0$ was arbitrary, $v$ minimizes $I$ on $\mathcal{S}_1^0$.

    To identify the support of this minimizer, choose $x_0 \in V$ such that
    \begin{equation*}
        \mu(x_0) \mathrm{e}^{2a(x_0)/b} = \min_{x \in V} \mu(x) \mathrm{e}^{2a(x)/b},
    \end{equation*}
    and define
    \begin{equation*}
        w_0(x) \coloneqq \begin{cases}
            \mu(x_0)^{-1/2}, & x = x_0,    \\[1mm]
            0,               & x \neq x_0.
        \end{cases}
    \end{equation*}
    Then $w_0 \in \mathcal{S}_1^0$, and
    \begin{equation*}
        I(w_0) = \frac{b}{4} \log\left(\mu(x_0) \mathrm{e}^{2a(x_0)/b}\right).
    \end{equation*}
    By \eqref{eq:finite-graph-limit-energy}, the inequality $I(v) \leq I(w_0)$ and $b > 0$ imply
    \begin{equation*}
        \sigma_K \leq \mu(x_0) \mathrm{e}^{2a(x_0)/b} = \min_{x \in V} \mu(x) \mathrm{e}^{2a(x)/b}.
    \end{equation*}
    Since $K$ is nonempty and each term in the sum defining $\sigma_K$ is at least this strictly positive minimum, $K$ must consist of a single vertex $x_*$. The same inequality then gives
    \begin{equation*}
        \mu(x_*) \mathrm{e}^{2a(x_*)/b} = \min_{x \in V} \mu(x) \mathrm{e}^{2a(x)/b}.
    \end{equation*}
    The mass constraint yields $v(x_*) = \mu(x_*)^{-1/2}$, and $v(x) = 0$ for every $x \neq x_*$.
\end{proof}

\begin{proof}[Theorem~\ref{thm:finite-graph-mass-limits}(ii)]
    Let $\{m_n\}$ be an arbitrary sequence of positive masses such that $m_n \to +\infty$, and define $v_n$, $\tau_n$ and $\theta_n$ as in part~(i). The identities \eqref{eq:finite-graph-mass-rescaled} and \eqref{eq:finite-graph-mass-multiplier} remain valid. The same compactness argument yields, up to a subsequence, $v_n \to v$ in $\mathcal{H}_0$ and uniformly on $V$, where $v \geq 0$ on $V$ and $\int_V v^2 \, \mathrm{d}\mu = 1$.

    Since $p > 2$, we have $\tau_n \to +\infty$. Using the same continuity argument as in part~(i), we obtain
    \begin{equation*}
        \frac{\theta_n}{\tau_n} = -\frac{2I(v_n)}{\tau_n} - \int_V |\nabla v_n|^p \, \mathrm{d}\mu \to -\int_V |\nabla v|^p \, \mathrm{d}\mu.
    \end{equation*}
    Using the convergence arguments in part~(i), we divide \eqref{eq:finite-graph-mass-rescaled} by $\tau_n$ and pass to the limit to obtain
    \begin{equation}
        -\Delta_p v(x) = \left(\int_V |\nabla v|^p \, \mathrm{d}\mu\right) v(x), \qquad x \in V.
        \label{eq:finite-graph-large-mass-limit}
    \end{equation}

    By the definition of $\Delta_p$, we have
    \begin{equation*}
        \begin{aligned}
            \int_V \Delta_p v \, \mathrm{d}\mu & = \frac{1}{2} \sum_{x \in V} \sum_{y \sim x} \omega_{xy} \left(|\nabla v|(x)^{p - 2} + |\nabla v|(y)^{p - 2}\right) (v(y) - v(x)) \\
                                               & = \frac{1}{2} \sum_{x \in V} \sum_{y \sim x} \omega_{yx} \left(|\nabla v|(y)^{p - 2} + |\nabla v|(x)^{p - 2}\right) (v(x) - v(y)) \\
                                               & = -\int_V \Delta_p v \, \mathrm{d}\mu.
        \end{aligned}
    \end{equation*}
    Hence $\int_V \Delta_p v \, \mathrm{d}\mu = 0$. Integrating \eqref{eq:finite-graph-large-mass-limit} over $V$, we obtain
    \begin{equation*}
        \left(\int_V |\nabla v|^p \, \mathrm{d}\mu\right) \left(\int_V v \, \mathrm{d}\mu\right) = 0.
    \end{equation*}
    Since $v \geq 0$ and $\int_V v^2 \, \mathrm{d}\mu = 1$, we have $\int_V v \, \mathrm{d}\mu > 0$. Hence $\int_{V}|\nabla v|^p \, \mathrm{d}\mu = 0$ which means $v$ is constant on $V$, and hence
    \begin{equation*}
        v(x) = \mu(V)^{-1/2}, \qquad x \in V.
    \end{equation*}
    By compactness of $\mathcal{S}_1^0$, the whole sequence $v_n$ therefore converges to $\mu(V)^{-1/2}$ in $\mathcal{H}_0$ and uniformly on $V$. As the sequence $\{m_n\}$ is arbitrary, this proves the asserted convergence as $m \to +\infty$.

    Integrate \eqref{eq:finite-graph-mass-rescaled} over $V$. Again using $\int_V \Delta_p v_n \, \mathrm{d}\mu = 0$ and $\int_V v_n \, \mathrm{d}\mu > 0$, we obtain
    \begin{equation*}
        \theta_n = \frac{\displaystyle b \int_V v_n \log v_n \, \mathrm{d}\mu - \int_V a(x) v_n \, \mathrm{d}\mu}{\displaystyle \int_V v_n \, \mathrm{d}\mu}.
    \end{equation*}
    Passing to the limit and using $v_n \to \mu(V)^{-1/2}$ uniformly on the finite graph, we find
    \begin{equation*}
        \theta_n \to -\frac{1}{\mu(V)} \int_V a(x) \, \mathrm{d}\mu - \frac{b}{2} \log \mu(V).
    \end{equation*}
    The definition of $\theta_n$ and the arbitrariness of $\{m_n\}$ give the stated limit of $\lambda_m - (b/2) \log m$ for the whole family.
\end{proof}

\section{The infinite graph case}
\label{sec:infinite-graph}

In this section, we consider infinite graphs with the graph and parameter conditions specified in the introduction. We use the spaces $\mathcal{W}$ and $\mathcal{H}$ and the energy functional $J$ defined there. For each $m > 0$, let
\begin{equation*}
    \mathcal{S}_m \coloneqq \left\{u \in \mathcal{H} \colon \int_V u^2 \, \mathrm{d}\mu = m\right\}.
\end{equation*}

\begin{lemma}
    \label{lem:finite-gradient-approximation}
    Assume (A1). Let $\mathcal{X}$ be either $\mathcal{W}$ or $\mathcal{H}$, and assume (A3) when $\mathcal{X} = \mathcal{H}$. For every $u \colon V \to \mathbb{R}$, we have $u \in \mathcal{X}$ if and only if $\|u\|_{\mathcal{X}} < +\infty$. In this case, there exists a sequence $\{\phi_n\} \subset C_c(V)$ such that
    \begin{equation*}
        \|\phi_n - u\|_{\mathcal{X}} \to 0.
    \end{equation*}
\end{lemma}

\begin{proof}
    Suppose $u \in \mathcal{X}$. By definition, there exists a sequence $\{\phi_n\} \subset C_c(V)$ such that
    \begin{equation*}
        \|\phi_n - u\|_{\mathcal{X}} \to 0.
    \end{equation*}
    The triangle inequality gives
    \begin{equation*}
        \|u\|_{\mathcal{X}} \leq \|\phi_n\|_{\mathcal{X}} + \|\phi_n - u\|_{\mathcal{X}} < +\infty.
    \end{equation*}

    Conversely, suppose $\|u\|_{\mathcal{X}} < +\infty$. Then $u \in L^2(V)$. Following \cite[Lemma~3.16(b)]{KellerLenzWojciechowski2021}, we define
    \begin{equation*}
        u^\varepsilon(x) \coloneqq
        \begin{cases}
            u(x) + \varepsilon, & u(x) < -\varepsilon, \\
            0, & |u(x)| \leq \varepsilon, \\
            u(x) - \varepsilon, & u(x) > \varepsilon,
        \end{cases}
    \end{equation*}
    where $\varepsilon > 0$. Since $\operatorname{supp}u^\varepsilon = \{x \in V \colon |u(x)| > \varepsilon\}$, by (A1) we have
    \begin{equation*}
        \sum_{x \in \operatorname{supp}u^\varepsilon} \mu_0\varepsilon^2 \leq \int_{\operatorname{supp}u^\varepsilon} u^2 \, \mathrm{d}\mu \leq \int_V u^2 \, \mathrm{d}\mu < +\infty.
    \end{equation*}
    Thus $u^\varepsilon \in C_c(V)$. By the definition of $u^\varepsilon$, we have
    \begin{equation*}
        |u - u^\varepsilon| \leq |u|, \qquad |\nabla u - \nabla u^\varepsilon| \leq |\nabla u|.
    \end{equation*}
    Since
    \begin{equation*}
        u - u^\varepsilon \to 0, \qquad \nabla u - \nabla u^\varepsilon \to 0 \quad \text{pointwise on } V
    \end{equation*}
    as $\varepsilon \to 0^+$, Lebesgue's dominated convergence theorem yields
    \begin{equation*}
        \|u^\varepsilon - u\|_{\mathcal{X}} \to 0.
    \end{equation*}
\end{proof}

\begin{lemma}
    \label{lem:H-reflexivity}
    Assume (A1). Let $\mathcal{X}$ be either $\mathcal{W}$ or $\mathcal{H}$, and assume (A3) when $\mathcal{X} = \mathcal{H}$. Then $\mathcal{X}$ is reflexive, and every bounded sequence $\{u_n\} \subset \mathcal{X}$ has a subsequence, still denoted by $\{u_n\}$, such that
    \begin{equation*}
        u_n \rightharpoonup u \quad \text{in } \mathcal{X}, \qquad u_n(x) \to u(x) \quad \text{for every } x \in V,
    \end{equation*}
    where $u \in \mathcal{X}$.
\end{lemma}

\begin{proof}
    Reflexivity follows as in the proof of \cite[Theorem~1.3(ii)]{Shao2025}.

    By reflexivity, every bounded sequence $\{u_n\} \subset \mathcal{X}$ has a subsequence, still denoted by $\{u_n\}$, such that $u_n \rightharpoonup u$ in $\mathcal{X}$ for some $u \in \mathcal{X}$.

    For each $x \in V$, take $\phi_x \coloneqq \mathbf{1}_{\{x\}} \in L^2(V)$. By the Cauchy--Schwarz inequality,
    \begin{equation*}
        \left|\int_V v \phi_x \, \mathrm{d}\mu\right| \leq \|v\|_{L^2(V)} \|\phi_x\|_{L^2(V)} = \sqrt{\mu(x)} \|v\|_{L^2(V)}, \qquad v \in \mathcal{X},
    \end{equation*}
    so the linear functional $v \mapsto \int_V v \phi_x \, \mathrm{d}\mu$ belongs to $\mathcal{X}^*$.

    Since $u_n \rightharpoonup u$ in $\mathcal{X}$, we have
    \begin{equation*}
        \mu(x) |u_n(x) - u(x)| = \left|\int_V (u_n - u) \phi_x \, \mathrm{d}\mu\right| \to 0.
    \end{equation*}
\end{proof}

\begin{lemma}
    \label{lem:fixed-mass-compactness}
    Assume (A1)--(A3). Let $\{u_n\} \subset \mathcal{S}_m$ satisfy
    \begin{equation*}
        \sup_n J(u_n) < +\infty.
    \end{equation*}
    Then there exist $u \in \mathcal{S}_m$ and a subsequence, still denoted by $\{u_n\}$, such that, for both $\mathcal{X} = \mathcal{W}$ and $\mathcal{X} = \mathcal{H}$,
    \begin{equation*}
        u_n \rightharpoonup u \quad \text{in } \mathcal{X}, \qquad u_n \to u \quad \text{in } L^2(V).
    \end{equation*}
\end{lemma}

\begin{proof}
    Set
    \begin{equation*}
        M \coloneqq \sqrt{\frac{m}{\mu_0}}.
    \end{equation*}
    The mass constraint and (A1) give
    \begin{equation*}
        |u_n(x)| \leq \sqrt{\frac{m}{\mu(x)}} \leq M.
    \end{equation*}
    There exists a constant $C > 0$, independent of $n$, such that
    \begin{equation*}
        \frac{1}{p} \int_V |\nabla u_n|^p \, \mathrm{d}\mu + \frac{1}{2} \int_V a(x) u_n^2 \, \mathrm{d}\mu \leq J(u_n) + \frac{bm}{2} \log M \leq C.
    \end{equation*}
    This estimate shows that $\{u_n\}$ is bounded in $\mathcal{H}$.

    By Lemma~\ref{lem:H-reflexivity}, there exist $u \in \mathcal{H}$ and a subsequence, still denoted by $\{u_n\}$, such that
    \begin{equation*}
        u_n \rightharpoonup u \quad \text{in } \mathcal{H}, \qquad u_n(x) \to u(x) \quad \text{for every } x \in V.
    \end{equation*}
    Since $\mathcal{H} \hookrightarrow \mathcal{W}$ continuously, the same subsequence also converges weakly to $u$ in $\mathcal{W}$.

    To prove strong convergence in $L^2(V)$, we follow the argument in the proof of \cite[Lemma~2.2]{GrigoryanLinYang2017}. By Fatou's lemma,
    \begin{equation*}
        \int_V a(x) u^2 \, \mathrm{d}\mu \leq \liminf_{n \to \infty} \int_V a(x) u_n^2 \, \mathrm{d}\mu \leq 2C.
    \end{equation*}
    Fix $R > 0$. By (A2), the set $F_R \coloneqq \{x \in V \colon a(x) \leq R\}$ is finite. Pointwise convergence therefore yields
    \begin{equation*}
        \int_{F_R} |u_n - u|^2 \, \mathrm{d}\mu = \sum_{x \in F_R} \mu(x) |u_n(x) - u(x)|^2 \to 0.
    \end{equation*}
    On $V \setminus F_R$, we have $a(x) > R$, so
    \begin{equation*}
        \int_{V \setminus F_R} |u_n - u|^2 \, \mathrm{d}\mu \leq \frac{2}{R} \int_V a(x) (u_n^2 + u^2) \, \mathrm{d}\mu \leq \frac{8C}{R}.
    \end{equation*}
    Combining these two estimates and letting $n \to \infty$ with $R$ fixed, we obtain
    \begin{equation*}
        \limsup_{n \to \infty} \|u_n - u\|_{L^2(V)}^2 \leq \frac{8C}{R}.
    \end{equation*}
    Letting $R \to +\infty$, we conclude that $u_n \to u$ in $L^2(V)$. Hence
    \begin{equation*}
        \int_V u^2 \, \mathrm{d}\mu = \lim_{n \to \infty} \int_V u_n^2 \, \mathrm{d}\mu = m,
    \end{equation*}
    so $u \in \mathcal{S}_m$.
\end{proof}

\begin{lemma}
    \label{lem:weak-test-extension}
    Assume (A1) and (A3), and let $(u, \lambda_m)$ be a solution of problem \eqref{eq:normalized-logarithmic-p-laplacian}. Then, for every $\phi \in \mathcal{H}$,
    \begin{equation*}
        \int_V |u \log u| |\phi| \, \mathrm{d}\mu < +\infty,
    \end{equation*}
    and
    \begin{equation*}
        \int_V |\nabla u|^{p - 2} \nabla u \cdot \nabla \phi \, \mathrm{d}\mu + \int_V [a(x) + \lambda_m] u \phi \, \mathrm{d}\mu = b \int_V u \phi \log u \, \mathrm{d}\mu.
    \end{equation*}
\end{lemma}

\begin{proof}
    We first prove that the logarithmic integral is absolutely convergent for every $\phi \in \mathcal{H}$.

    By (A1), (A3), and the mass constraint,
    \begin{equation*}
        \|u\|_{L^\infty(V)} \leq \sqrt{\frac{m}{\mu_0}}, \qquad \|\phi\|_{L^2(V)} \leq a_0^{-1/2} \|\phi\|_{\mathcal{H}} \quad \text{for every } \phi \in \mathcal{H}.
    \end{equation*}
    H\"older's inequality shows that $Lu$, defined by
    \begin{equation*}
        \langle Lu, \phi \rangle \coloneqq \int_V |\nabla u|^{p - 2} \nabla u \cdot \nabla \phi \, \mathrm{d}\mu + \int_V [a(x) + \lambda_m] u \phi \, \mathrm{d}\mu,
    \end{equation*}
    belongs to $\mathcal{H}^*$. For $\eta \in C_c(V)$, we have $\bigl\||\eta|\bigr\|_{\mathcal{H}} \leq \|\eta\|_{\mathcal{H}}$. Testing \eqref{eq:infinite-graph-variational-identity} with $|\eta|$ gives
    \begin{equation*}
        b \int_V |u \log u| |\eta| \, \mathrm{d}\mu = 2b \int_V u (\log u)_+ |\eta| \, \mathrm{d}\mu - \langle Lu, |\eta| \rangle \leq C \|\eta\|_{\mathcal{H}},
    \end{equation*}
    where $(\log u)_+ \coloneqq \max\{\log u, 0\}$ and $C$ is independent of $\eta$.

    For $\phi \in \mathcal{H}$, choose $\phi_n \in C_c(V)$ with $\phi_n \to \phi$ in $\mathcal{H}$, hence pointwise on $V$. Fatou's lemma gives
    \begin{equation*}
        b \int_V |u \log u| |\phi| \, \mathrm{d}\mu \leq \liminf_{n \to \infty} b \int_V |u \log u| |\phi_n| \, \mathrm{d}\mu \leq C \|\phi\|_{\mathcal{H}} < +\infty.
    \end{equation*}
    Hence \eqref{eq:infinite-graph-variational-identity} holds for every $\phi \in \mathcal{H}$.
\end{proof}

We now prove Theorem~\ref{thm:infinite-graph-existence} by minimizing $J$ on $\mathcal{S}_m$.

\begin{proof}[Theorem~\ref{thm:infinite-graph-existence}]
    For any $x_0 \in V$, taking $\phi_m \coloneqq \sqrt{m/\mu(x_0)}\,\mathbf{1}_{\{x_0\}} \in C_c(V) \cap \mathcal{S}_m$ gives $\inf\limits_{u \in \mathcal{S}_m} J(u) \leq J(\phi_m) < +\infty$. 
    
    For every $u \in \mathcal{S}_m$, (A1) gives $|u(x)| \leq \sqrt{m/\mu_0}$ on $V$, so (A3) and \eqref{eq:energy-functional} yield
    \begin{equation*}
        J(u) \geq -\frac{b}{2} \int_V u^2 \log |u| \, \mathrm{d}\mu \geq -\frac{bm}{4} \log\frac{m}{\mu_0} > -\infty.
    \end{equation*}
    Choose a minimizing sequence $\{u_n\} \subset \mathcal{S}_m$ such that
    \begin{equation*}
        J(u_n) \to \inf_{u \in \mathcal{S}_m} J(u).
    \end{equation*}
    By Lemma~\ref{lem:fixed-mass-compactness}, after passing to a subsequence, there exists $u_m \in \mathcal{H}$ such that $u_n \to u_m$ in $L^2(V)$. Consequently,
    \begin{equation*}
        \int_V u_m^2 \, \mathrm{d}\mu = \lim_{n \to \infty} \int_V u_n^2 \, \mathrm{d}\mu = m,
    \end{equation*}
    so $u_m \in \mathcal{S}_m$. Moreover,
    \begin{equation*}
        \inf_{u \in \mathcal{S}_m} J(u) \leq J(u_m) \leq \liminf_{n \to \infty} J(u_n) = \inf_{u \in \mathcal{S}_m} J(u).
    \end{equation*}
    Hence $J(u_m) = \inf\limits_{u \in \mathcal{S}_m} J(u)$.

    The reverse triangle inequality gives $|\nabla(|u_m|)| \leq |\nabla u_m|$, so Lemma~\ref{lem:finite-gradient-approximation} yields $|u_m| \in \mathcal{S}_m$. As in the proof of Theorem~\ref{thm:finite-graph-existence}, $J(|u_m|) \leq J(u_m)$, so $|u_m|$ is also a minimizer. Replacing $u_m$ by $|u_m|$ if necessary, we may assume that $u_m \geq 0$ on $V$.

    Since $J(u_m) < +\infty$, the defining formula for $I$ in Section~\ref{sec:finite-graph} also defines a finite value $I(u_m)$. We set
    \begin{equation*}
        \lambda_m \coloneqq -\frac{1}{m} \left(\int_V |\nabla u_m|^p \, \mathrm{d}\mu + 2I(u_m)\right) \in \mathbb{R}.
    \end{equation*}
    For any $\phi \in C_c(V)$, we use the perturbation $su_m + t\phi$ as in the proof of \cite[Lemma~3.5]{Shao2026a}. Let $g(s,t) \coloneqq J(su_m + t\phi)$ and
    \begin{equation*}
        F(s,t) \coloneqq \int_V (su_m + t\phi)^2 \, \mathrm{d}\mu - m.
    \end{equation*}
    For $s > 0$,
    \begin{equation*}
        J(su_m) = \frac{s^p}{p} \int_V |\nabla u_m|^p \, \mathrm{d}\mu + s^2 I(u_m) + \frac{bm}{4} s^2 (1 - 2\log s).
    \end{equation*}
    Since $\phi$ has finite support and the graph is locally finite, $g(s,t) - J(su_m)$ is a finite sum of $C^1$ functions, so $g$ is $C^1$ near $(1,0)$. Direct differentiation gives
    \begin{equation}
        \frac{\partial F}{\partial s}(1,0) = 2m > 0, \qquad \frac{\partial g}{\partial s}(1,0) = \int_V |\nabla u_m|^p \, \mathrm{d}\mu + 2I(u_m) = -m\lambda_m.
        \label{eq:infinite-graph-s-derivatives}
    \end{equation}

    By the minimality of $u_m$, the point $(1,0)$ minimizes $g$ subject to $F = 0$. By \eqref{eq:infinite-graph-s-derivatives}, the Lagrange multiplier method gives $c \in \mathbb{R}$ such that
    \begin{equation*}
        \frac{\partial g}{\partial s}(1,0) + c\frac{\partial F}{\partial s}(1,0) = 0, \qquad \frac{\partial g}{\partial t}(1,0) + c\frac{\partial F}{\partial t}(1,0) = 0.
    \end{equation*}
    Substituting \eqref{eq:infinite-graph-s-derivatives} into the first identity gives $c = \lambda_m/2$. And the second identity yields
    \begin{equation}
        \int_V |\nabla u_m|^{p - 2} \nabla u_m \cdot \nabla \phi \, \mathrm{d}\mu + \int_V [a(x) + \lambda_m] u_m \phi \, \mathrm{d}\mu = b \int_V u_m \phi \log u_m \, \mathrm{d}\mu.
        \label{eq:infinite-graph-variational-equation}
    \end{equation}

    Taking $\phi$ to be the indicator function of each singleton in \eqref{eq:infinite-graph-variational-equation} and using integration by parts, we obtain
    \begin{equation*}
        -\Delta_p u_m(x) + [a(x) + \lambda_m] u_m(x) = b u_m(x) \log u_m(x), \qquad x \in V.
    \end{equation*}
    By the same argument as in the proof of Theorem~\ref{thm:finite-graph-existence}, we have $u_m > 0$ on $V$, and $(u_m, \lambda_m)$ is a solution of \eqref{eq:normalized-logarithmic-p-laplacian}.
\end{proof}

We now turn to the asymptotic behavior of solutions as the mass tends to zero or infinity. We first prove Theorem~\ref{thm:infinite-graph-mass-limits}(i), concerning the limit $m \to 0^+$.

\begin{proof}[Theorem~\ref{thm:infinite-graph-mass-limits}(i)]
    Let $m_n \to 0^+$ and define $v_n$, $\tau_n$ and $\theta_n$ as in the proof of Theorem~\ref{thm:finite-graph-mass-limits}(i). Since $v_n \in \mathcal{S}_1$, (A1) gives
    \begin{equation}
        0 < v_n(x) \leq \mu(x)^{-1/2} \leq \mu_0^{-1/2}, \qquad x \in V.
        \label{eq:infinite-graph-small-mass-pointwise-bound}
    \end{equation}
    By Lemma~\ref{lem:weak-test-extension}, we have
    \begin{equation*}
        \int_V v_n^2 |\log v_n| \, \mathrm{d}\mu < +\infty,
    \end{equation*}
    so \eqref{eq:finite-graph-mass-rescaled} and \eqref{eq:finite-graph-mass-multiplier} remain valid. By \eqref{eq:infinite-graph-small-mass-pointwise-bound}, we have $\log v_n(x) \leq -(1/2) \log\mu_0$, and hence
    \begin{equation}
        2I(v_n) \geq \int_V a(x) v_n^2 \, \mathrm{d}\mu + \frac{b}{2} \log\mu_0.
        \label{eq:infinite-graph-unit-mass-I-bound}
    \end{equation}
    By (A3), \eqref{eq:finite-graph-mass-multiplier} and \eqref{eq:infinite-graph-unit-mass-I-bound}, $\{\theta_n\}$ is bounded above. Thus, along a subsequence, either $\theta_n \to \theta \in \mathbb{R}$ or $\theta_n \to -\infty$. Since $\|v_n\|_{L^2(V)} = 1$, we may also assume that
    \begin{equation*}
        v_n \rightharpoonup v \quad \text{in } L^2(V), \qquad v_n(x) \to v(x) \geq 0 \quad \text{for every } x \in V.
    \end{equation*}

    Assume $\theta_n \to -\infty$ and suppose, for contradiction, that $v(x) > 0$ for some $x \in V$. Since $\tau_n \Delta_p v_n(x) \to 0$ by local finiteness, \eqref{eq:finite-graph-mass-rescaled} gives
    \begin{equation*}
        b v(x) \log v(x) - a(x) v(x) = \lim_{n \to \infty} \theta_n v_n(x) = -\infty,
    \end{equation*}
    a contradiction. Hence $v \equiv 0$.

    Suppose now that $\theta_n \to \theta \in \mathbb{R}$. By \eqref{eq:finite-graph-mass-multiplier} and \eqref{eq:infinite-graph-unit-mass-I-bound},
    \begin{equation*}
        \int_V a(x) v_n^2 \, \mathrm{d}\mu \leq -\theta_n - \frac{b}{2} \log\mu_0 \leq C
    \end{equation*}
    for some $C > 0$ independent of $n$. For $R > 0$, the set $F_R \coloneqq \{x \in V \colon a(x) \leq R\}$ is finite by (A2). Since $a(x) > R$ outside $F_R$, we obtain
    \begin{equation*}
        \int_{F_R} v^2 \, \mathrm{d}\mu = \lim_{n \to \infty} \int_{F_R} v_n^2 \, \mathrm{d}\mu \geq 1 - \frac{C}{R}.
    \end{equation*}
    Letting $R \to +\infty$ gives $\|v\|_{L^2(V)} \geq 1$. Weak convergence gives $\|v\|_{L^2(V)} \leq 1$, so $v_n \to v$ in $L^2(V)$ and $\|v\|_{L^2(V)} = 1$.

    Let $K \coloneqq \operatorname{supp} v$. Passing to the limit in \eqref{eq:finite-graph-mass-rescaled} at each $x \in K$ and using (A3), we obtain
    \begin{equation*}
        v(x) = \mathrm{e}^{(a(x) + \theta)/b} \geq \mathrm{e}^{(a_0 + \theta)/b} > 0, \qquad x \in K.
    \end{equation*}
    Since $\|v\|_{L^2(V)} = 1$ and (A1) holds, $K$ is nonempty and finite. The mass constraint \eqref{eq:small-mass-normalization} then gives \eqref{eq:small-mass-explicit-limit}.

    Suppose now that each $u_m$ is a minimizer obtained in the proof of Theorem~\ref{thm:infinite-graph-existence}. By (A1) and (A2), for every $R \in \mathbb{R}$ the set
    \begin{equation*}
        \bigl\{x \in V \colon \mu(x) \mathrm{e}^{2a(x)/b} \leq R\bigr\}
    \end{equation*}
    is finite. Hence $\mu(x) \mathrm{e}^{2a(x)/b}$ attains its minimum at some $x_0 \in V$. Set $w_0 \coloneqq \mu(x_0)^{-1/2}\mathbf{1}_{\{x_0\}}$.

    By \eqref{eq:infinite-graph-small-mass-pointwise-bound}, we have
    \begin{equation*}
        2I(v_n) \geq a(x_0) + \frac{b}{2} \log \mu(x_0) = 2I(w_0).
    \end{equation*}
    As in the proof of Theorem~\ref{thm:finite-graph-mass-limits}(i), we obtain \eqref{eq:finite-graph-minimizer-comparison} with $w = w_0$, namely
    \begin{equation*}
        I(v_n) + \frac{\tau_n}{p} \int_V |\nabla v_n|^p \, \mathrm{d}\mu \leq I(w_0) + \frac{\tau_n}{p} \int_V |\nabla w_0|^p \, \mathrm{d}\mu.
    \end{equation*}
    Hence $I(v_n) \to I(w_0)$ as $\tau_n \to 0$. By \eqref{eq:finite-graph-mass-multiplier}, $\theta_n \to -2I(w_0)$. Comparing this limit with \eqref{eq:small-mass-explicit-limit} yields
    \begin{equation*}
        \sigma_K = \mu(x_0) \mathrm{e}^{2a(x_0)/b} = \min_{x \in V} \mu(x) \mathrm{e}^{2a(x)/b}.
    \end{equation*}
    Every summand in $\sigma_K$ is at least this positive minimum. Hence $K$ consists of a single minimizing vertex $x_*$, and \eqref{eq:small-mass-minimizer-limit} follows.
\end{proof}

For the large-mass limit in Theorem~\ref{thm:infinite-graph-mass-limits}(ii), we use the following consequence of \eqref{eq:infinite-graph-gradient-infimum}. By homogeneity and the density of $C_c(V)$ in $\mathcal{W}$,
\begin{equation}
    \int_V |\nabla u|^p \, \mathrm{d}\mu \geq \kappa \|u\|_{L^2(V)}^p, \qquad u \in \mathcal{W}.
    \label{eq:infinite-graph-gradient-lower-bound}
\end{equation}

We also need the following inequality.

\begin{lemma}[Discrete Picone inequality {\cite[Proposition~4.2]{BrascoFranzina2014}}]
    \label{lem:discrete-picone}
    Let $v > 0$ and $\phi \geq 0$ on $V$. For every edge $xy \in E$,
    \begin{equation*}
        (v(y) - v(x)) \left(\frac{\phi(y)^2}{v(y)} - \frac{\phi(x)^2}{v(x)}\right) \leq (\phi(y) - \phi(x))^2.
    \end{equation*}
\end{lemma}

\begin{proof}[Theorem~\ref{thm:infinite-graph-mass-limits}(ii)]
    Define $v_n$, $\tau_n$ and $\theta_n$ as in the proof of Theorem~\ref{thm:finite-graph-mass-limits}(i). Then $\tau_n \to +\infty$, $v_n \in \mathcal{S}_1$, and \eqref{eq:infinite-graph-small-mass-pointwise-bound} holds. By Lemma~\ref{lem:weak-test-extension}, the integrals defining $I(v_n)$ are finite, and \eqref{eq:finite-graph-mass-rescaled} and \eqref{eq:finite-graph-mass-multiplier} remain valid. Hence, we obtain
    \begin{equation}
        -\Delta_p v_n + \frac{a(x)}{\tau_n} v_n + \frac{\theta_n}{\tau_n} v_n = \frac{b}{\tau_n} v_n \log v_n,
        \label{eq:infinite-graph-large-mass-rescaled}
    \end{equation}
    \begin{equation}
        \frac{\theta_n}{\tau_n} = -\int_V |\nabla v_n|^p \, \mathrm{d}\mu - \frac{2I(v_n)}{\tau_n}.
        \label{eq:infinite-graph-large-mass-multiplier}
    \end{equation}
    Also, \eqref{eq:infinite-graph-unit-mass-I-bound} gives
    \begin{equation}
        2I(v_n) \geq a_0 + \frac{b}{2} \log\mu_0.
        \label{eq:infinite-graph-large-mass-I-lower-bound}
    \end{equation}

    Since $\|v_n\|_{L^2(V)} = 1$, we may assume up to a subsequence that $v_n \rightharpoonup v$ in $L^2(V)$, with $\|v\|_{L^2(V)} \leq 1$. Testing against $\phi_x = \mathbf{1}_{\{x\}}$ as in the proof of Lemma~\ref{lem:H-reflexivity} gives $v_n(x) \to v(x) \geq 0$ for every $x \in V$.

    Suppose $v \not\equiv 0$ and fix $x_0$ with $v(x_0) > 0$. By local finiteness and \eqref{eq:infinite-graph-large-mass-rescaled},
    \begin{equation*}
        \frac{\theta_n}{\tau_n} = \frac{\Delta_p v_n(x_0)}{v_n(x_0)} + \frac{b}{\tau_n} \log v_n(x_0) - \frac{a(x_0)}{\tau_n} \to \theta, \qquad \theta \coloneqq \frac{\Delta_p v(x_0)}{v(x_0)}.
    \end{equation*}
    Passing to the limit at each vertex in \eqref{eq:infinite-graph-large-mass-rescaled} gives
    \begin{equation}
        -\Delta_p v + \theta v = 0 \quad \text{on } V.
        \label{eq:infinite-graph-large-mass-theta-limit}
    \end{equation}
    As in the proof of Theorem~\ref{thm:finite-graph-existence}, we obtain $v > 0$. Since $\|v_n\|_{L^2(V)} = 1$, it follows from \eqref{eq:infinite-graph-gradient-lower-bound}, \eqref{eq:infinite-graph-large-mass-multiplier} and \eqref{eq:infinite-graph-large-mass-I-lower-bound} that
    \begin{equation}
        \kappa \leq \|\nabla v_n\|_{L^p(V)}^p \leq -\frac{\theta_n + a_0 + (b/2) \log\mu_0}{\tau_n} \to -\theta.
        \label{eq:infinite-graph-large-mass-gradient-bound}
    \end{equation}
    Thus $\kappa \leq -\theta$. By density we may test \eqref{eq:infinite-graph-large-mass-theta-limit} with $v$, obtaining
    \begin{equation}
        \|\nabla v\|_{L^p(V)}^p = -\theta \|v\|_{L^2(V)}^2 > 0,
        \label{eq:infinite-graph-gradient-testing}
    \end{equation}
    since (A1) gives $\mu(V) = +\infty$ and thus $v$ is nonconstant.

    For $\phi \in C_c(V)$ with $\|\phi\|_{L^2(V)} = 1$, we test \eqref{eq:infinite-graph-large-mass-theta-limit} with $\phi^2/v$, following Brasco and Franzina~\cite{BrascoFranzina2014}. Applying Lemma~\ref{lem:discrete-picone} to $v$ and $|\phi|$, and then H\"older's inequality, we obtain
    \begin{equation}
        -\theta \leq \int_V |\nabla v|^{p - 2} |\nabla\phi|^2 \, \mathrm{d}\mu \leq \|\nabla v\|_{L^p(V)}^{p - 2} \|\nabla\phi\|_{L^p(V)}^2.
        \label{eq:infinite-graph-large-mass-test-bound}
    \end{equation}

    Taking the infimum over $\phi$ in \eqref{eq:infinite-graph-large-mass-test-bound} and using \eqref{eq:infinite-graph-large-mass-gradient-bound}, \eqref{eq:infinite-graph-gradient-testing} and $\|v\|_{L^2(V)} \leq 1$, we obtain
    \begin{equation*}
        -\theta \leq \kappa \|v\|_{L^2(V)}^{p - 2} \leq \kappa \leq -\theta.
    \end{equation*}
    Since $-\theta > 0$ and $p > 2$, equality gives $\kappa = -\theta$ and $\|v\|_{L^2(V)} = 1$. Thus $v$ satisfies \eqref{eq:infinite-graph-large-mass-limit} and $\theta_n/\tau_n \to -\kappa$.

    Since $v_n \rightharpoonup v$ in $L^2(V)$ and $\|v_n\|_{L^2(V)} = \|v\|_{L^2(V)} = 1$, we have $v_n \to v$ in $L^2(V)$.

    Finally, suppose that each $u_m$ is a minimizer obtained in the proof of Theorem~\ref{thm:infinite-graph-existence} and that $\kappa > 0$. Fix $w \in C_c(V)$ with $\|w\|_{L^2(V)} = 1$. By \eqref{eq:finite-graph-minimizer-comparison} and \eqref{eq:infinite-graph-large-mass-I-lower-bound},
    \begin{equation*}
        \|\nabla v_n\|_{L^p(V)}^p \leq \|\nabla w\|_{L^p(V)}^p + \frac{C_w}{\tau_n},
    \end{equation*}
    where $C_w$ is independent of $n$.

    Suppose, for contradiction, that $v \equiv 0$. Then $v_n \rightharpoonup 0$ in $L^2(V)$, so $\|v_n + w\|_{L^2(V)}^2 \to 2$. Since $\nabla w$ has finite support and $\nabla v_n \to 0$ pointwise,
    \begin{equation*}
        \|\nabla(v_n + w)\|_{L^p(V)}^p - \|\nabla v_n\|_{L^p(V)}^p \to \|\nabla w\|_{L^p(V)}^p.
    \end{equation*}
    Applying \eqref{eq:infinite-graph-gradient-lower-bound} to $v_n + w$ and letting $n \to \infty$, we obtain
    \begin{equation*}
        2^{p/2}\kappa \leq 2\|\nabla w\|_{L^p(V)}^p.
    \end{equation*}
    Taking the infimum over $w$ gives $2^{p/2}\kappa \leq 2\kappa$, contradicting $p > 2$ and $\kappa > 0$. Thus $v \not\equiv 0$.
\end{proof}

%% ── Acknowledgements ─────────────────────────────────
% \input{sections/acknowledgements}

%% ── Bibliography ──────────────────────────────────────
\bibliographystyle{abbrv}
\bibliography{refs}
\end{document}